\documentclass[12pt]{article}

\usepackage[center]{titlesec}
\usepackage[leqno]{amsmath}

\usepackage{amsfonts}

\begin{document}

\def\no{\noindent\bf}

\title{{\bf{\small{JORDAN HIGHER ALL-DERIVABLE POINTS IN NEST ALGEBRAS}}\footnotetext{Communicated by Nannan Zhen.}}}

\author{\small{Nannan Zhen and Jun Zhu}}

\date{}
\normalsize
\maketitle
\begin{quote}
 {\bf Abstract.}{ Let $\mathcal{N}$ be a non-trivial and complete nest on a Hilbert space
 $H$. Suppose $d=\{d_n: n\in N\}$ is a group of linear mappings from Alg$\mathcal{N}$ into itself. We say that $d=\{d_n: n\in N\}$ is a Jordan higher derivable mapping at a given point $G$ if $d_{n}(ST+ST)=\sum\limits_{i+j=n}\{d_{i}(S)d_{j}(T)+d_{j}(T)d_{i}(S)\}$ for any $S,T\in Alg \mathcal{N}$ with $ST=G$. An element $G\in Alg \mathcal{N}$ is called a Jordan higher all-derivable point if every Jordan higher derivable mapping at $G$ is a higher derivation. In this paper, we mainly prove that any given point $G$ of Alg$\mathcal{N}$ is a Jordan higher all-derivable point. This extends some results in \cite{Chen11} to the case of higher derivations.}\\
\end{quote}
\normalsize
\section{I{\small{NTRODUCTION}} {\small{AND}} P{\small{RELIMINARIES}}}
~~
Let $\mathcal{A}$ be an algebra. A linear mapping $\delta$ from $\mathcal{A}$ into itself is called a derivation if $\delta(ST)=\delta(S)T+S\delta(T)$ for any $S,T\in\mathcal{A}$ and is said to be a Jordan derivation if $\delta(ST+TS)=\delta(S)T+S\delta(T)+\delta(T)S+T\delta(S)$ for any $S,T\in\mathcal{A}$. We say that a mapping $\delta$ is Jordan derivable at a given point $G\in \mathcal{A}$ if $\delta(ST)=\delta(S)T+S\delta(T)+\delta(T)S+T\delta(S)$ and $G$ is called a Jordan all-derivable point of $\mathcal{A}$ if every Jordan derivable mapping at $G$ is a derivation. Suppose that $d=\{d_n: n\in N\}$ is a group of linear mappings from $\mathcal{A}$ into itself and $d_{0}$ is the identical mapping. We say that $d=\{d_n: n\in N\}$ is Jordan higher derivable at a given point $G$ if $d_{n}(ST+ST)=\sum\limits_{i+j=n}\{d_{i}(S)d_{j}(T)+d_{j}(T)d_{i}(S)\}$ for any $S,T\in\mathcal{A}$ with $ST=G$, and $G$ is called a Jordan higher all-derivable point of $\mathcal{A}$ if every Jordan higher derivable mapping at $G$ is a higher derivation, that is $d_{n}(ST)=\sum\limits_{i+j=n}\{d_{i}(S)d_{j}(T)\}$ for any $S,T\in\mathcal{A}$.

With the development of derivation, higher derivation has attracted much attention of mathematicians as an active subject of research in algebras. Generally speaking, there are two directions in the study of the local actions of derivations of operator algebras. One is the well known local derivation problem. The other is to study conditions under which derivations of operator algebras can be completely determined by the action on some sets of operators. It is obvious that a linear map is a higher derivation if and only if it is higher derivable at all points. It is natural and interesting to ask the question whether or not a linear map is a derivation if it is Jordan higher derivable only at one given point.

We describe some of the results related to ours. In \cite{Chen11} ,Chen proved that any $G\neq0$ is a Jordan all-derivable point in nest algebras. Zhao S and Zhu J pointed that $G=0$ is a Jordan all-derivable point in nest algebras in \cite{Zhao10}. Jing, Lu and Li \cite{Jing02} showed that every derivable mapping $\varphi$ at $0$ with $\varphi(I)=0$ on nest algebras is a derivation. In \cite{Xiao10}, Z. Xiao and F. Wei gave the proof of the fact that any Jordan higher derivation on a triangular algebra is a higher derivation.

In this paper, $\mathcal{N}$ is a non-trivial and complete nest on a Hilbert space $H$. $Alg\mathcal{N}=\{A\in B(H):AP\subseteq P, \forall P\in \mathcal{N}\}$ is an algebra. Let $N\in \mathcal{N}$ with $0\subset N\subset H$. Then we can get the orthogonal decomposition $H=N\oplus N^{\bot}$. In this way, we can write $G=\left [\begin{array}{cc}
  D & E \\
  0 & F \\
\end{array}\right ]$, where $D\in Alg\mathcal{N}_{N}$, $E\in B(N^{\bot},N)$ and $F\in Alg\mathcal{N}_{N^{\bot}}$ $(\mathcal{N}_{N}=\{M\bigcap N:M\in N)\},\mathcal{N}_{N^{\bot}}=\{M\bigcap N^{\bot}:M\in N\})$. All the identitical mappings in the proof are represented by $I$ and $\lambda$ is a positive real number for convenient writing.

\section{J{\small{ORDAN}} H{\small{IGHER}} A{\small{LL-DERIVABLE}} P{\small{OINTS}} {\small{IN}} N{\small{EST}} A{\small{LGEBRAS}}}
~~
In this section, we assume that $d=\{d_n: n\in N\}$ is a Jordan higher derivable linear mapping at $G$ from $Alg\mathcal{N}$ into itself. We only need to prove that $d=\{d_n: n\in N\}$ is a higher derivation.\\
~

{\bf{Theorem 2.1}} {\it Let $\mathcal{N}$ be a non-trivial and complete nest on a Hilbert space $H$. Any element of $Alg\mathcal{N}$ is a Jordan higher all-derivable point}.

{\bf Proof.} For any $X\in Alg\mathcal{N}_{N}$, $Y\in B(N^{\bot},N)$, $Z\in Alg\mathcal{N}_{N^{\bot}}$, we can write $$d_{n}(\left [\begin{array}{cc}
  X & Y \\
  0 & Z \\
\end{array}\right ])=\left [\begin{array}{cc}
  A_{n1}(X)+B_{n1}(Y)+C_{n1}(Z) & A_{n2}(X)+B_{n2}(Y)+C_{n2}(Z)  \\
  0 & A_{n3}(X)+B_{n3}(Y)+C_{n3}(Z)  \\
\end{array}\right ]$$, Where $A_{ij}$,$B_{ij}$ and $C_{ij}$ are linear mappings on $Alg\mathcal{N}_{N}$ , $B(N^{\bot},N) $ and $Alg\mathcal{N}_{N}^{\bot}$ ,respectively. It is clear that $A_{01}(X)=X$, $A_{02}(X)=0$, $A_{03}(X)=0$, $B_{01}(Y)=0$, $B_{02}(Y)=Y$, $B_{03}(Y)=0$, $C_{01}(Z)=0$, $C_{02}(Z)=0$ and $C_{03}(Z)=Z$.

Let $S=\left [\begin{array}{cc}
  X & Y \\
  0 & Z \\
\end{array}\right ]$, $T=\left [\begin{array}{cc}
  U & V \\
  0 & W \\
\end{array}\right ]$ for any $X,U\in Alg\mathcal{N}_{N}$, $Y,V\in B(N^{\bot},N)$, $Z,W\in Alg\mathcal{N}_{{N}^{\bot}}$ with $XU=D$, $XV+YW=E$ and $ZW=F$, then $ST=G$. So we have \\$d_{n}(ST+TS)=\\
\left [\begin{array}{ccc}
  A_{n1}(D)+A_{n1}(UX)+B_{n1}(E)+ &     A_{n2}(D)+A_{n2}(UX)+B_{n2}(E)+\\
  B_{n1}(UY+VZ)+C_{n1}(F)+C_{n1}(WZ)& B_{n2}(UY+VZ)+C_{n2}(F)+C_{n2}(WZ)\\
0 & A_{n3}(D)+A_{n3}(UX)+B_{n3}(E)+ \\~& B_{n3}(UY+VZ)+C_{n3}(F)+C_{n3}(WZ)\\
\end{array}\right ]\\
=\sum\limits_{i+j=n}\{d_{i}(S)d_{j}(T)+d_{j}(T)d_{i}(S)\}=\sum\limits_{i+j=n}\\
\{\left [\begin{array}{ccc}
  A_{i1}(X)+B_{i1}(Y)&A_{i2}(X)+B_{i2}(Y)\\
  +C_{i1}(Z) & +C_{i2}(Z)\\
  0 & A_{i3}(X)+B_{i3}(Y)\\~& +C_{i3}(Z)\\
\end{array}\right ]\bullet \left [\begin{array}{ccc}
  A_{j1}(U)+B_{j1}(V)&A_{j2}(U)+B_{j2}(V)\\
  +C_{j1}(W) & +C_{j2}(W)\\
  0 & A_{j3}(U)+B_{j3}(V)\\~& +C_{j3}(W)\\
\end{array}\right ]\\
+ \left [\begin{array}{ccc}
  A_{j1}(U)+B_{j1}(V)&A_{j2}(U)+B_{j2}(V)\\
  +C_{j1}(W) & +C_{j2}(W)\\
  0 & A_{j3}(U)+B_{j3}(V)\\~& +C_{j3}(W)\\
\end{array}\right ]\bullet \left [\begin{array}{ccc}
  A_{i1}(X)+B_{i1}(Y)&A_{i2}(X)+B_{i2}(Y)\\
  +C_{i1}(Z) & +C_{i2}(Z)\\
  0 & A_{i3}(X)+B_{i3}(Y)\\~& +C_{i3}(Z)\
\end{array}\right ] \}$. This implies that
\begin{equation}
A_{n1}(D)+A_{n1}(UX)+B_{n1}(E)+B_{n1}(UY+VZ)+C_{n1}(F)+C_{n1}(WZ)\\
\end{equation}
$~~~~=\sum\limits_{i+j=n}\{A_{i1}(X)A_{j1}(U)+A_{i1}(X)B_{j1}(V)+A_{i1}(X)C_{j1}(W)+B_{i1}(Y)A_{j1}(U)\\
 ~~~~~~~+B_{i1}(Y)B_{j1}(V)+B_{i1}(Y)C_{j1}(W)+C_{i1}(Z)A_{j1}(U)+C_{i1}(Z)B_{j1}(V)\\
 ~~~~~~~+C_{i1}(V)C_{j1}(W)+A_{j1}(U)A_{i1}(X)+A_{j1}(U)B_{i1}(Y)+A_{j1}(U)C_{i1}(Z)\\
 ~~~~~~~+B_{j1}(V)A_{i1}(X)+B_{j1}(V)B_{i1}(Y)+B_{j1}(V)C_{i1}(Z)+C_{j1}(W)A_{i1}(X)\\
 ~~~~~~~+C_{j1}(W)B_{i1}(Y)+C_{j1}(W)C_{i1}(Z)\}$,
\begin{equation}
A_{n2}(D)+A_{n2}(UX)+B_{n2}(E)+B_{n2}(UY+VZ)+C_{n2}(F)+C_{n2}(WZ)
\end{equation}
$ ~~~~=\sum\limits_{i+j=n}\{A_{i1}(X)A_{j2}(U)+A_{i1}(X)B_{j2}(V)+A_{i1}(X)C_{j2}(W)+B_{i1}(Y)A_{j2}(U)\\
 ~~~~~~~+B_{i1}(Y)B_{j2}(V)+B_{i1}(Y)C_{j2}(W)+C_{i1}(Z)A_{j2}(U)+C_{i1}(Z)B_{j2}(V)\\
 ~~~~~~~+C_{i1}(Z)C_{j2}(W)+A_{i2}(X)A_{j3}(U)+A_{i2}(X)B_{j3}(V)+A_{i2}(X)C_{j3}(W)\\
 ~~~~~~~+B_{i2}(Y)A_{j3}(U)+B_{i2}(Y)B_{j3}(V)+B_{i2}(Y)C_{j3}(W)+C_{i2}(Z)A_{j3}(U)\\
 ~~~~~~~+C_{i2}(Z)B_{j3}(V)+C_{i2}(Z)C_{j3}(W)+A_{j1}(U)A_{i2}(X)+A_{j1}(U)B_{i2}(Y)\\
 ~~~~~~~+A_{j1}(U)C_{i2}(Z)+B_{j1}(V)A_{i2}(X)+B_{j1}(V)B_{i2}(Y)+B_{j1}(V)C_{i2}(Z)\\
 ~~~~~~~+C_{j1}(W)A_{i2}(X)+C_{j1}(W)B_{i2}(Y)+C_{j1}(W)C_{i2}(Z)+A_{j2}(U)A_{i3}(X)\\
 ~~~~~~~+A_{j2}(U)B_{i3}(Y)+A_{j2}(U)C_{i3}(Z)+B_{j2}(V)A_{i3}(X)+B_{j2}(V)B_{i3}(Y)\\
 ~~~~~~~+B_{j2}(V)C_{i3}(Z)+C_{j2}(W)A_{i3}(X)+C_{j2}(W)B_{i3}(Y)+C_{j2}(W)C_{i3}(Z)\}$\\and
\begin{equation}
A_{n3}(D)+A_{n3}(UX)+B_{n3}(E)+B_{n3}(UY+VZ)+C_{n3}(F)+C_{n3}(WZ)
\end{equation}
$ ~~~~=\sum\limits_{i+j=n}\{A_{i3}(X)A_{j3}(U)+A_{i3}(X)B_{j3}(V)+A_{i3}(X)C_{j3}(W)+B_{i3}(Y)A_{j3}(U)\\
 ~~~~~~~+B_{i3}(Y)B_{j3}(V)+B_{i3}(Y)C_{j3}(W)+C_{i3}(Z)A_{j3}(U)+C_{i3}(Z)B_{j3}(V)\\
 ~~~~~~~+C_{i3}(Z)C_{j3}(W)+A_{j3}(U)A_{i3}(X)+A_{j3}(U)B_{i3}(Y)+A_{j3}(U)C_{i3}(Z)\\
 ~~~~~~~+B_{j3}(V)A_{i3}(X)+B_{j3}(V)B_{i3}(Y)+B_{j3}(V)C_{i3}(Z)+C_{j3}(W)A_{i3}(X)\\
 ~~~~~~~+C_{j3}(W)B_{i3}(Y)+C_{j3}(W)C_{i3}(Z)\}$.

 {\bf Case1.} $G\neq0$.

 {\bf Step1.} We show that $G_{n1}(W)=0$ for any $W\in Alg\mathcal{N}_{N^{\bot}}$.

 Taking $X=\lambda^{-1}X, Y=\lambda Y, Z=\lambda Z, U=\lambda U, V=\lambda V$ and $W=\lambda^{-1} W$  with $XU=D$, $XV+YW=E$ and $ZW=F$ in Eq.(1), it follows that
\begin{equation}
A_{n1}(D)+A_{n1}(UX)+B_{n1}(E)+\lambda^{2}B_{n1}(UY+VZ)+C_{n1}(F)+C_{n1}(WZ)
\end{equation}
$~~~=\sum\limits_{i+j=n}\{A_{i1}(X)A_{j1}(U)+A_{i1}(X)B_{j1}(V)+\frac{1}{\lambda^{2}}A_{i1}(X)C_{j1}(W)+\lambda^{2}B_{i1}(Y)A_{j1}(U)\\
~~~~~~~~+\lambda^{2}B_{i1}(Y)B_{j1}(V)+B_{i1}(Y)C_{j1}(W)+\lambda^{2}C_{i1}(Z)A_{j1}(U)+\lambda^{2}C_{i1}(Z)B_{j1}(V)\\
~~~~~~~~+C_{i1}(Z)C_{j1}(W)+A_{j1}(U)A_{i1}(X)+\lambda^{2}A_{j1}(U)B_{i1}(Y)+\lambda^{2}A_{j1}(U)C_{i1}(Z)\\
~~~~~~~~+B_{j1}(V)A_{i1}(X)+\lambda^{2}B_{j1}(V)B_{i1}(Y)+\lambda^{2}B_{j1}(V)C_{i1}(Z)+\frac{1}{\lambda^{2}}C_{j1}(W)A_{i1}(X)\\
~~~~~~~~+C_{j1}(W)B_{i1}(Y)+C_{j1}(W)C_{i1}(Z)\}$.\\
Multiplying Eq.(4) by $\lambda^{2}$ and let $\lambda \rightarrow 0$, then $\sum\limits_{i+j=n}\{A_{i1}(X)C_{j1}(W)+C_{j1}(W)A_{i1}(X)\}$\\
 =0. It is clearly established when $n=0$. When $n=1$, we can get $XC_{11}(W)+C_{11}(W)X=0$ for any $X\in Alg\mathcal{N}$ with $XU=D$. Taking $X=I$, then $C_{11}(W)=0$. We assume that $C_{m1}(W)=0$ for all $0\leq m<n$. In fact, after simplifying the equation, we have $XC_{n1}(W)+C_{n1}(W)X=0$. Taking $X=I$, then $C_{n1}(W)=0$. Thus $C_{n1}(W)=0$ for any $W\in Alg\mathcal{N}_{N^{\bot}}$.

{\bf Step2.} We show that $A_{n1}(XU+UX)=\sum\limits_{i+j=n}\{A_{i1}(X)A_{j1}(U)+A_{j1}(U)A_{i1}(X)\}$ for any $X,U\in Alg\mathcal{N}_{N}$ with $XU=D$ and $B_{n1}(V)=0$ for any $V\in B(N^{\bot},N)$.

Taking $X=\lambda^{-1}X, Y=\lambda E, Z=\lambda F, U=\lambda U, V=0$ and $W=\lambda^{-1}I$ with $XU=D$ in Eq.(1), then
\begin{equation}
A_{n1}(D)+A_{n1}(UX)+B_{n1}(E)+\lambda^{2}B_{n1}(UE)\nonumber
\end{equation}
~~~~~$=\sum\limits_{i+j=n}\{A_{i1}(X)A_{j1}(U)+\lambda^{2}B_{i1}(E)A_{j1}(U)+A_{j1}(U)A_{i1}(X)+\lambda^{2}A_{j1}(U)B_{i1}(E)\}$. \\ Dividing the above equation by $\lambda^{2}$ and $\lambda\rightarrow\infty$, then $B_{n1}(UE)=\sum\limits_{i+j=n}\{B_{i1}(E)A_{j1}(U)+A_{j1}(U)B_{i1}(E)\}$. Taking $X=D$ and $U=I$, thus $B_{n1}(E)=\sum\limits_{i+j=n}\{B_{i1}(E)A_{j1}(I)+A_{j1}(I)B_{i1}(E)\}$. By induction, we get $B_{n1}(E)=0$. It follows that $A_{n1}(XU+UX)=\sum\limits_{i+j=n}\{A_{i1}(X)A_{j1}(U)+A_{j1}(U)A_{i1}(X)\}$ with $XU=D$. Now the simplified Eq. (4) is $\sum\limits_{i+j=n}\{A_{i1}(X)B_{j1}(V)+B_{j1}(V)A_{i1}(X)\}=0$. Applying mathematical induction, we gain that $B_{n1}(V)=0$ for any $V\in B(N^{\bot},N)$.

{\bf Step3.} We show that  $A_{n3}(X)=0$ for any $X\in Alg\mathcal{N}_{N}$ and $B_{n3}(Y)=0$ for any $Y\in B(N^{\bot},N)$.

For any $Y\in B(N^{\bot},N)$, taking $X=\lambda I, Y=Y, Z=-\lambda^{-1}F, U=\lambda^{-1}D, V=Y+\lambda^{-1}E$ and $W=-\lambda I$ in Eq.(3), then
\begin{equation}
A_{n3}(2D)+B_{n3}(E)+\frac{1}{\lambda}B_{n3}(DY-YF)-\frac{1}{\lambda^{2}}B_{n3}(EF)+C_{n3}(2F)\nonumber
\end{equation}
$~~~~~~~=\sum\limits_{i+j=n}\{A_{i3}(I)A_{j3}(D)+{\lambda}A_{i3}(I)B_{j3}(Y)+A_{i3}(I)B_{j3}(E)-{\lambda^{2}}A_{i3}(I)C_{j3}(I)\\
~~~~~~~~~~~~+\frac{1}{\lambda}B_{i3}(Y)A_{j3}(D)+B_{i3}(Y)B_{j3}(Y)+\frac{1}{\lambda}B_{i3}(Y)B_{j3}(E)-{\lambda}B_{i3}(Y)C_{j3}(I)\\
~~~~~~~~~~~~-\frac{1}{\lambda^{2}}C_{i3}(F)A_{j3}(D)-\frac{1}{\lambda}C_{i3}(F)B_{j3}(Y)-\frac{1}{\lambda^{2}}C_{i3}(F)B_{j3}(E)+C_{i3}(F)C_{j3}(I)\\
~~~~~~~~~~~~+A_{j3}(D)A_{i3}(I)+{\lambda}B_{j3}(Y)A_{i3}(I)+B_{j3}(E)A_{i3}(I)-\lambda^{2}C_{j3}(I)A_{i3}(I)\\
~~~~~~~~~~~~+\frac{1}{\lambda}A_{j3}(D)B_{j3}(Y)+B_{j3}(Y)B_{i3}(Y)+\frac{1}{\lambda}B_{j3}(E)B_{i3}(Y)-{\lambda}C_{j3}(I)B_{i3}(Y)\\
~~~~~~~~~~~~-\frac{1}{\lambda^{2}}A_{j3}(D)C_{i3}(F)-\frac{1}{\lambda}B_{j3}(Y)C_{i3}(F)-\frac{1}{\lambda^{2}}B_{j3}(E)C_{i3}(F)+C_{j3}(I)C_{i3}(F)\}$.\\
With the randomicity of $\lambda$, we infer that
\begin{equation}
A_{n3}(D)+A_{n3}(D)+B_{n3}(E)+C_{n3}(2F)\nonumber
\end{equation}
$~~~~~~~=\sum\limits_{i+j=n}\{A_{i3}(I)A_{j3}(D)+A_{i3}(I)B_{j3}(E)+B_{i3}(Y)B_{j3}(Y)+C_{i3}(F)C_{j3}(I)\\
~~~~~~~~~~+A_{j3}(D)A_{i3}(I)+B_{j3}(E)A_{i3}(I)+B_{j3}(Y)B_{i3}(Y)+C_{j3}(I)C_{i3}(F)\}$ \\for any $Y\in B(N^{\bot},N)$. The two sides will be, in fact, equal when $Y=0$. So $\sum\limits_{i+j=n}\{B_{i3}(Y)B_{j3}(Y)+B_{j3}(Y)B_{i3}(Y)\}=0$. Applying mathematical induction, we claim $B_{n3}(Y)=0$ for any $Y\in B(N^{\bot},N)$.

Taking $X=\lambda^{-1}X, Y=\lambda E, Z=\lambda F, U=\lambda U, V=0$ and $W=\lambda^{-1}I$  with $XU=D$ in Eq.(3),then
\begin{equation}
A_{n3}(D)+A_{n3}(UX)+C_{n3}(F)+C_{n3}(F)\nonumber
\end{equation}
$~~~~=\sum\limits_{i+j=n}\{A_{i3}(X)A_{j3}(U)+\frac{1}{\lambda^{2}}A_{i3}(X)C_{j3}(I)+\lambda^{2}C_{i3}(F)A_{j3}(U)+C_{i3}(F)C_{j3}(I)\\
~~~~~~~~~~+A_{j3}(U)A_{i3}(X)+\lambda^{2}A_{j3}(U)C_{i3}(F)+\frac{1}{\lambda^{2}}C_{j3}(I)A_{i3}(X)+C_{j3}(I)C_{i3}(F)\}$. \\ Multiplying the above equation by $\lambda^{2}$ and let $\lambda\rightarrow0$,then$\sum\limits_{i+j=n}\{A_{i3}(X)C_{j3}(I)+C_{j3}(I)A_{i3}(X)\}=0$. Similarly available $A_{n3}(X)=0$ for any $X\in Alg\mathcal{N}_{N}$.

{\bf Step4.} We show that $A_{n1}(I)=0$ and $C_{n3}(I)=0$ for all $n>0$.

{\bf Case\textcircled{1}.} $D$ and $F$, at least one is not $0$.

Taking $X=-X, Y=Y, Z=Z, U=-U, V=-V$ and $W=W$ with $XU=D$, $XV+YW=E$ and $ZW=F$ in Eq.(2), we get a new Eq.(2).The new one together with the original one yield
\begin{equation}
A_{n2}(D)+A_{n2}(UX)+B_{n2}(E)+C_{n2}(F)+C_{n2}(WZ)
\end{equation}
$ ~~~~~~~~~~~~~~~~=\sum\limits_{i+j=n}\{A_{i1}(X)A_{j2}(U)+A_{i1}(X)B_{j2}(V)+B_{i2}(Y)C_{j3}(W)\\
~~~~~~~~~~~~~~~~~~~~+C_{i2}(Z)C_{j3}(W)+A_{j1}(U)A_{i2}(X)+C_{j2}(W)C_{i3}(Z)\}$.\\
For any $V\in B(N^{\bot},N)$, taking $X=I, Y=E-V, Z=F, U=D, V=V$ and $W=I$ in Eq.(5), then we can get $\sum\limits_{i+j=n}\{A_{i1}(I)B_{j2}(V)-B_{i2}(V)C_{j3}(I)\}=0$. Also by mathematical induction, we prove that
\begin{equation}
A_{n1}(I)V=VC_{n3}(I).
\end{equation}
We can easily get $A_{n1}(2D)=\sum\limits_{i+j=n}\{A_{i1}(D)A_{j1}(I)+A_{j1}(I)A_{i1}(D)\}$ from step 2 and $C_{n3}(2F)=\sum\limits_{i+j=n}\{C_{i3}(F)C_{j3}(I)+C_{j3}(I)C_{i3}(F)\}$ by taking $Z=F, W=I$ in Eq.(3). Then, using induction, $DA_{n1}(I)+A_{n1}(I)D=0$ and $FC_{n3}(I)+C_{n3}(I)F=0$ are obtained by assuming $C_{m3}(I)=0, A_{m1}(I)=0, 0<m<n$.

By $A_{11}(I)=0$ and $C_{13}(I)=0$ from [1,p.56], we get $FC_{23}(I)+C_{23}(I)F=0, DA_{21}(I)+A_{21}(I)D=0$. Combining with (6) and the condition of $D$ and $F$, we have $A_{21}(I)=0,C_{23}(I)=0$ by the same method in [1,p.56]. Also, with the help of induction, we point that $A_{n1}(I)=0, C_{23}(I)=0$ for all $n>0$.

{\bf Case\textcircled{2}.} $D=0, F=0, E\neq0$

Taking $X=0, Y=E, Z=0, U=I, V=V$ and $W=I$ in Eq.(2), then we see\\
\begin{equation}
B_{n2}(E)+B_{n2}(E)=\sum\limits_{i+j=n}\{B_{i2}(E)C_{j3}(I)+A_{j1}(I)B_{i2}(E)\}.\nonumber
\end{equation}
 By $A_{11}(I)=0$ and $C_{13}(I)=0$ from [1,p.56], the above equation implies $EC_{23}(I)+A_{21}(I)E=0$ when $n=2$. We assume $A_{m1}(I)=0, C_{m3}(I)=0$ for all $0<m<n$. In fact, after simplifying, we get $EC_{n3}(I)+A_{n1}(I)E=0$.

By the same way in [1,p.56], we can prove that $A_{n1}(I)=0, C_{n3}(I)=0$ for all $n>0$.

{\bf Step5.} We show that $A_{n2}(X)=-\sum\limits_{i+j=n}A_{i1}(X)C_{j2}(I)$ for any $X\in Alg\mathcal{N}_{N}$ and $C_{n2}(W)=-\sum\limits_{i+j=n}A_{i2}(I)C_{j3}(W)$ for any $W\in Alg\mathcal{N}_{N^{\bot}}$.

For any invertible $X$, taking $X=X, Y=\lambda E, Z=\lambda F, U=X^{-1}D, V=0$ and $W=\lambda^{-1}I$ in Eq.(2), it follows that \\
\begin{equation}
A_{n2}(X^{-1}DX+D)+B_{n2}(E)+\lambda B_{n2}(X^{-1}DE)+C_{n2}(2F)\nonumber
\end{equation}
~~~~~~~~~~~~~~$=\sum\limits_{i+j=n}\{A_{i1}(X)A_{j2}(X^{-1}D)+\frac{1}{\lambda}A_{i1}(X)C_{j2}(I)+\frac{1}{\lambda}A_{i2}(X)C_{j3}(I)\\
~~~~~~~~~~~~~~+{\lambda}A_{j1}(X^{-1}D)C_{i2}(F)+{\lambda}A_{j2}(X^{-1}D)C_{i3}(F)+{\lambda}A_{j1}(X^{-1}D)B_{i2}(E)\\
~~~~~~~~~~~~~~+A_{j1}(X^{-1}D)A_{i2}(X)+B_{i2}(E)C_{j3}(I)+C_{i2}(F)C_{j3}(I)+C_{j2}(I)C_{i3}(F)\}$.\\
 Multiplying the above equation by $\lambda$ and let $\lambda\rightarrow0$, then $\sum\limits_{i+j=n}\{A_{i1}(X)C_{j2}(I)+A_{i2}(X)C_{j3}(I)\}=0$. In fact, from step 4, we can get $A_{n2}(X)=-\sum\limits_{i+j=n}A_{i1}(X)C_{j2}(I)$. So, for any invertible operator $X\in Alg\mathcal{N}_{N}$, $A_{n2}(X)=-\sum\limits_{i+j=n}A_{i1}(X)C_{j2}(I)$.  We
can get the fact that $A_{n2}(X)=-\sum\limits_{i+j=n}A_{i1}(X)C_{j2}(I)$ for any $X\in Alg\mathcal{N}_{N}$ from [1,Lemma4.1].

For any invertible $W\in Alg\mathcal{N}_{N^{\bot}}$, taking $X=I, Y=0, Z=\lambda^{-1}FW^{-1}, U=D, V=E$ and $W=\lambda W$ in Eq.(2), that is
\begin{equation}
A_{n2}(2D)+B_{n2}(E)+\frac{1}{\lambda}B_{n2}(EFW^{-1})+C_{n2}(F)+C_{n2}(WFW^{-1})\nonumber
\end{equation}
$~~~~~=\sum\limits_{i+j=n}\{A_{i1}(I)A_{j2}(D)+A_{i1}(I)B_{j2}(E)+\lambda A_{i1}(I)C_{j2}(W)+\lambda A_{i2}(I)C_{j3}(W)\\
~~~~~~~~~+C_{i2}(FW^{-1})C_{j3}(W)+A_{j1}(D)A_{i2}(I)+\frac{1}{\lambda}A_{j1}(D)C_{i2}(FW^{-1})\\
~~~~~~~~~+\frac{1}{\lambda}A_{j2}(D)C_{i3}(FW^{-1})+\frac{1}{\lambda}B_{j2}(E)C_{i3}(FW^{-1})+C_{j2}(W)C_{i3}(FW^{-1})\}$.\\
 Dividing the above equation by $\lambda$ and let $\lambda\rightarrow\infty$, then $\sum\limits_{i+j=n}\{A_{i1}(I)C_{j2}(W)+A_{i2}(I)C_{j3}(W)\}=0$. By step 4, we have $C_{n2}(W)=-\sum\limits_{i+j=n}A_{i2}(I)C_{j3}(W)$ for any invertible $W\in Alg\mathcal{N}_{N^{\bot}}$. By [1,Lemma4.1], for any $W\in Alg\mathcal{N}_{N^{\bot}}$, $C_{n2}(W)=-\sum\limits_{i+j=n}A_{i2}(I)C_{j3}(W)$.

{\bf Step6.} We show that $C_{n2}(ZW+WZ)=\sum\limits_{i+j=n}\{C_{i2}(Z)C_{j3}(W)+C_{j2}(W)C_{i3}(Z)\}$ with $ZW=F$ and $\sum\limits_{i+j=n}\{A_{i1}(I)C_{j2}(W)+A_{i2}(I)C_{j3}(W)\}=0$ for any $X\in Alg\mathcal{N}_{N}$ and $W\in Alg\mathcal{N}_{N^{\bot}}$.

We have $\sum\limits_{i+j=n}\{A_{i1}(I)C_{j2}(W)+A_{i2}(I)C_{j3}(W)\}=0$ by step 5.

Taking $X=I, Y=0, Z=Z, U=D, V=E$ and $W=W$ with $ZW=F$ in Eq.(5), then
\begin{equation}
A_{n2}(D)+A_{n2}(D)+B_{n2}(E)+C_{n2}(F)+C_{n2}(WZ)\nonumber
\end{equation}
$~~~~~~~~~~~~~~~~~~~=\sum\limits_{i+j=n}\{A_{i1}(I)A_{j2}(D)+A_{i1}(I)B_{j2}(E)\\
~~~~~~~~~~~~~~~~~~~~+C_{i2}(Z)C_{j3}(W)+A_{j1}(D)A_{i2}(I)+C_{j2}(W)C_{i3}(Z)\}$.\\
 We can easily get
$C_{n2}(ZW+WZ)=\sum\limits_{i+j=n}\{C_{i2}(Z)C_{j3}(W)+C_{j2}(W)C_{i3}(Z)\}$ with $ZW=F$ after simplifying the above equation.

{\bf Step7.} We show that both $A_{n1}(\cdot)$ and $C_{n3}(\cdot)$ are derivations.

For any invertible $W\in Alg\mathcal{N}_{N^{\bot}}$, taking $X=I, Y=Y, Z=FW^{-1}, U=D, V=E-YW$ and $W=W$ in Eq.(5), then we get
\begin{equation}
B_{n2}(YW)=\sum\limits_{i+j=n}\{B_{i2}(Y)C_{j3}(W)\}.\nonumber
\end{equation}
By [1,Lemma4.1], we know that $B_{n2}(YW)=\sum\limits_{i+j=n}\{B_{i2}(Y)C_{j3}(W)\}$ for any $W\in Alg\mathcal{N}_{N^{\bot}}$. For any $W_{1},W_{2}\in Alg\mathcal{N}_{N^{\bot}}$, we have $B_{n2}(YW_{1}W_{2})=\sum\limits_{i+j=n}\{B_{i2}(Y)C_{j3}(W_{1}W_{2})\}$. On the other hand, $B_{n2}(YW_{1}W_{2})=\sum\limits_{i+j=n}\{B_{i2}(YW_{1})C_{j3}(W_{2})\}=\sum\limits_{i+j+k=n}\{B_{i2}(Y)C_{k3}(W_{1})\\C_{j3}(W_{2})\}$. The above equations imply that $C_{n3}(W_{1}W_{2})=\sum\limits_{i+j=n}\{C_{i3}(W_{1})C_{j3}(W_{2})\}$.Hence $C_{i3}(\cdot)$ is a derivation.

For any invertible $X\in Alg\mathcal{N}_{N}$, taking $X=X, Y=E-XV, Z=F, U=X^{-1}D,V=V$ and $W=I$ in Eq.(5), it follows that $B_{n2}(XV)=\sum\limits_{i+j=n}\{A_{i1}(X)B_{j2}(V)\}$. Similarly, we can prove that $A_{n1}(\cdot)$ is a derivation.

{\bf Step8.} We show that $d_{n}(\cdot)$ is a higher derivation.

For arbitrary $S=\left [\begin{array}{cc}
  X & Y \\
  0 & Z \\
\end{array}\right ]$ and $T=\left [\begin{array}{cc}
  U & V \\
  0 & W \\
\end{array}\right ]$ in $Alg\mathcal{N}$, we only need to prove that $d_{n}(ST)=\sum\limits_{i+j=n}d_{i}(S)d_{j}(T)$. By the above steps, we calculate\\
$\sum\limits_{i+j=n}d_{i}(S)d_{j}(T)\\
=\sum\limits_{i+j=n}\{\left [\begin{array}{cc}
  A_{i1}(X)& A_{i2}(X)+B_{i2}(Y)\\
  &+C_{i2}(Z) \\
  0 & C_{i3}(Z)\\
\end{array}\right ]\bullet\left [\begin{array}{cc}
  A_{j1}(U) & A_{j2}(U)+B_{j2}(V)\\
  &+C_{j2}(W) \\
  0 & C_{j3}(W) \\
\end{array}\right ]\}\\
=\left [\begin{array}{cc}
  A_{n1}(XU)&\sum\limits_{i+j=n}\{A_{i1}(X)A_{j2}(U)+A_{i1}(X)B_{j2}(V)\\
  &+B_{i2}(Y)C_{j3}(W)+C_{i2}(Z)C_{j3}(W)\}\\
  0 & C_{n3}(ZW)\\
\end{array}\right ]\\
=\left [\begin{array}{cc}
  A_{n1}(XU)&B_{n2}(XV)+B_{n2}(YW)+\sum\limits_{i+j=n}\{A_{i1}(X)A_{j2}(U)+C_{i2}(Z)C_{j3}(W)\}\\
  0 & C_{n3}(ZW)\\
\end{array}\right ]\\
=\left [\begin{array}{cc}
  A_{n1}(XU)&B_{n2}(XV)+B_{n2}(YW)+\sum\limits_{i+j+K=n}\{-A_{i1}(X)A_{m1}(U)C_{k2}(I)\\
              &-A_{i2}(I)C_{m3}(Z)C_{j3}(W)\}\\
  0 & C_{n3}(ZW)\\
\end{array}\right ]\\
=\left [\begin{array}{cc}
  A_{i1}(XU)&B_{n2}(XV)+B_{n2}(YW)\\
  &+\sum\limits_{i+j+K=n}\{-A_{i1}(XU)C_{k2}(I)-A_{i2}(I)C_{m+j3}(ZW)\}\\
  0 & C_{n3}(ZW)\\
\end{array}\right ]\\
=\left [\begin{array}{cc}
  A_{i1}(XU)&B_{n2}(XV)+B_{n2}(YW)+A_{n2}(XU)+C_{n2}(ZW)\}\\
  0 & C_{n3}(ZW)\\
\end{array}\right ]=d_{n}(ST)$.\\

This completes the proof of case 1.

{\bf Case2.} $G=0$.

{\bf Step1.} We show that $B_{n1}(V)=0$ and $B_{n3}(V)=0$ for any $V\in B(N^{\bot},N)$.

Taking $X=0, Y=V, Z=0, U=X, V=0$ and $W=0$ in Eq.(1), Eq.(2) and Eq.(3) respectively, it follows that
\begin{equation}
B_{n1}(XV)=\sum\limits_{i+j=n}\{B_{i1}(V)A_{j1}(X)+A_{j1}(X)B_{i1}(V)\},
\end{equation}
\begin{eqnarray}
B_{n2}(XV)=\sum\limits_{i+j=n}\{B_{i1}(V)A_{j2}(X)+B_{i2}(V)A_{j3}(X)\nonumber\\
+A_{j2}(X)B_{i3}(V)+A_{j1}(X)B_{i2}(V)\},
\end{eqnarray}
\begin{equation}
B_{n3}(XV)=\sum\limits_{i+j=n}\{B_{i3}(V)A_{j3}(X)+A_{j3}(X)B_{i3}(V)\}.
\end{equation}
Taking $X=I$ in Eq.(7) and Eq.(9), by mathematical induction, we prove $B_{n1}(V)=0$ and $B_{n3}(V)=0$, respectively. The Eq.(8) can be simplified to
\begin{equation}
B_{n2}(XV)=\sum\limits_{i+j=n}\{B_{i2}(V)A_{j3}(X)+A_{j1}(X)B_{i2}(V)\}.
\end{equation}

{\bf Step2.} We show that $A_{n3}(X)=0$ for any $X\in Alg\mathcal{N}_{N}$ and $C_{n1}(W)=0$ for any $W\in Alg\mathcal{N}_{N^{\bot}}$.

Taking $X=0, Y=V, Z=W, U=X, V=0$ and $W=0$ in Eq.(1), Eq.(2) and Eq.(3) respectively, hence the following three equations hold
\begin{equation}
0=\sum\limits_{i+j=n}\{C_{i1}(W)A_{j1}(X)+A_{j1}(X)C_{i1}(W)\},
\end{equation}
\begin{eqnarray}
B_{n2}(XV)=\sum\limits_{i+j=n}\{C_{i1}(W)A_{j2}(X)+C_{i2}(W)A_{j3}(X)+B_{i2}(V)A_{j3}(X)\nonumber\\
+A_{j1}(X)C_{i2}(W)+A_{j1}(X)B_{i2}(V)+A_{j2}(X)C_{i3}(W)\},
\end{eqnarray}
\begin{equation}
0=\sum\limits_{i+j=n}\{C_{i3}(W)A_{j3}(X)+A_{j3}(X)C_{i3}(W)\}.
\end{equation}
Taking $X=I$ and $W=I$ in Eq.(11) and Eq.(13) respectively, applying mathematical induction , we can prove $C_{n1}(W)=0$ and $A_{n3}(X)=0$. Eq.(10) together with Eq.(12) yield
\begin{equation}
B_{n2}(XV)=\sum\limits_{i+j=n}A_{j1}(X)B_{i2}(V)
\end{equation}
and
\begin{equation}
\sum\limits_{i+j=n}\{A_{j1}(X)C_{i2}(W)+A_{j2}(X)C_{i3}(W)\}=0.
\end{equation}

{\bf Step3.} We show that $C_{n2}(W)=\sum\limits_{i+j=n}C_{i2}(I)C_{j3}(W)$ for any $W\in Alg\mathcal{N}_{N^{\bot}}$, $A_{n1}(\cdot)$ is is a derivation and $A_{n1}(I)=0$ for all $n>0$.

Similarly, from step 7 in case 1, we can prove $A_{n1}(\cdot)$ is a derivation by Eq.(14). Taking $X=I$ in Eq.(14), by mathematical induction, we get $A_{n1}(I)=0$ for all $n>0$. Taking $X=I$ in Eq.(15), then we have $C_{n2}(W)=\sum\limits_{i+j=n}C_{i2}(I)C_{j3}(W)$.

{\bf Step4.} We show that $A_{n2}(X)=-\sum\limits_{i+j=n}A_{i1}(X)C_{j2}(I)$ for any $X\in Alg\mathcal{N}_{N}$, $C_{n3}(\cdot)$ is a derivation and $C_{n3}(I)=0$ for all $n>0$.

Taking $X=\lambda I, Y=-V, Z=0, U=0, V=V$ and $W=\lambda I$ in Eq.(2), then
\begin{equation}
\sum\limits_{i+j=n}\{\lambda A_{i1}(I)B_{j2}(V)+\lambda^{2}A_{i1}(I)C_{j2}(I)+\lambda^{2}A_{i2}(I)C_{j3}(I)-\lambda B_{i2}(V)C_{j3}(I)\}=0.\nonumber
\end{equation}
Dividing the above equation by $\lambda$ and let $\lambda\rightarrow0$, then we get $\sum\limits_{i+j=n}\{A_{i1}(I)B_{j2}(V)-B_{i2}(V)C_{j3}(I)\}=0$. In fact, by step 3, $B_{n2}(V)=\sum\limits_{i+j=n}\{B_{i2}(V)C_{j3}(I)\}$. Using mathematical induction, we prove $C_{n3}(I)=0$ for all $n>0$. Taking $W=I$ in Eq.(15), that is $A_{n2}(X)=-\sum\limits_{i+j=n}A_{i1}(X)C_{j2}(I)$.

Taking $X=0, Y=-X^{-1}VW, Z=W, U=X, V=V$ and $W=0$ in Eq.(2), it follows that
\begin{equation}
\sum\limits_{i+j=n}\{-A_{j1}(X)B_{i2}(VW)+A_{j1}(X)C_{i2}(W)+A_{j2}(X)C_{i3}(W)+B_{j2}(V)C_{i3}(W)\}=0.\nonumber
\end{equation}
Taking $X=I$, with the help of Eq.(15), we get $B_{n2}(VW)=\sum\limits_{i+j=n}B_{j2}(V)C_{i3}(W)$. Similarly, we can prove that $C_{n3}(\cdot)$ is a derivation.

Similarly available $d=\{d_{n}:n\in N\}$ is a higher derivation.

\def\refname{\hfil {R\small{EFERENCES}}\hfil}

 ~\\{\small{Nannan Zhen and Jun Zhu\\
Institute of Mathematics,\\
Hangzhou Dianzi University, \\
Hangzhou 310018,\\
P.R. China\\
E-mail:nannanzhen66@sina.com}}

{\small

\begin{thebibliography}{99}


\bibitem{Chen11} Yunhe Chen, The depicting of mappings in operator algebra,
the doctoral dissertation of ECUST,2011.

\bibitem{Gong10} Ming Gong and Jun Zhu, Jordan multiplicative mappings at some points on matrix algebras,
Journal of Advanced Research in Pure Mathematics, 2010, 2(4): 84-93.

\bibitem{Jing02}Wu Jing, Shijie Lu and Pengting Li, Characterizations of derivations on some operator algebras,
Bull. Austral.Math. Soc. 66(2002),227-232.

\bibitem{Hou10} Xiaofeng Qi and Jinchuan Hou, Characterizations of derivations of Banach space nest algebras:
all-derivable point, Linear Algebra Appl. 432(2010),3183-3200.

\bibitem{Xiao10} Zhankui Xiao and Feng Wei, Jordan higher derivations of triangular algebras ,
Linear Algebra Appl,432(2010)2615-2622.

\bibitem{Zhao10} Sha Zhao and Jun Zhu, Jordan all-derivable points in the algebra of all upper triangular matrices[J],
 Linear Algebra Appl,2010,433(11-12)1922-1938.
 
\bibitem{Zhu09} Jun Zhu, Changping Xiong and Lin Zhang, All-derivable points in matrix algebras[J],
Linear Algebra Appl,2009,430(8-9):2070- 2079.







\end{thebibliography}
}
\end{document}